\input amstex
\magnification=1200
\documentstyle{amsppt}
\NoRunningHeads
\NoBlackBoxes
\define\BP{\bold P}
\define\Train{\operatorname{Train}}
\define\Mantle{\operatorname{Mantle}}
\define\Voile{\operatorname{Voile}}
\define\Trinion{\operatorname{Trinion}}
\define\Polynion{\operatorname{Polynion}}
\define\Vect{\operatorname{Vect}}
\define\CVect{\operatorname{\Bbb CVect}}
\define\Diff{\operatorname{Diff}}
\define\vir{\operatorname{vir}}
\define\Cvir{\operatorname{\Bbb Cvir}}
\define\Vir{\operatorname{Vir}}
\define\Gf{\operatorname{Gf}}
\define\PSLTWO{\operatorname{PSL}(2,\Bbb R)}
\define\PSLTWOC{\operatorname{PSL}(2,\Bbb C)}
\define\sltwo{\operatorname{\frak s\frak l}(2,\Bbb C)}
\define\sltwor{\operatorname{\frak s\frak l}(2,\Bbb R)}
\define\Hom{\operatorname{Hom}}
\define\End{\operatorname{End}}
\define\Mor{\operatorname{Mor}}
\define\SU{\operatorname{\bold S\bold U}}
\define\LDiff{\operatorname{\Cal LDIff}}
\define\LNer{\operatorname{\Cal LNer}}
\define\Ner{\operatorname{Ner}}
\define\sgn{\operatorname{sgn}}
\define\Rot{\operatorname{Rot}}
\define\Vrt{\operatorname{Vert}}
\define\PK{\operatorname{\Cal P\Cal K}}
\topmatter
\title On the infinite-dimensional hidden symmetries. II.
$q_R$-conformal modular functors.
\endtitle
\author Denis V. Juriev
\endauthor
\affil\eightpoint\rm
ul.Miklukho-Maklaya 20-180, Moscow 117437 Russia.\linebreak
E-mail: denis\@juriev.msk.ru\linebreak
\tenpoint\linebreak
January, 06, 1997\linebreak
E-version: funct-an/9701009
\endaffil
\abstract The article is devoted to the $q_R$-conformal modular functors,
which being ``deformations'' of the conformal modular functors (the projective
representations of the category $\Train(\Diff_+(S^1))$, the train of the
group $\Diff_+(S^1)$ of all orientation preserving diffeomorphisms of a
circle) in the class of all projective modular functors (the projective
representations of the category $\Train(\PSLTWO)$, the train of the projective
group $\PSLTWO$), may be regarded as their ``Berezin quantizations''.
\endabstract
\endtopmatter
\document
This paper being the continuation of the first part [1] belongs to the series
of articles supplemental to [2], and also lies in lines of the general ideology
exposed in the review [3]. The main purpose of the activity, which has its
origin and motivation presumably in the author's applied researches [4] on the
interactively controlled systems (i.e. the controlled systems, in which the
control is coupled with unknown or uncompletely known feedbacks), is to
explicate the essentially {\sl infinite-dimensional\/} aspects of the hidden
symmetries, which appear in the representation theory of the finite dimensional
Lie algebras and related algebraic structures. The relations between the
control and the representation theories will be discussed in [5]. The present
series is organized as a sequence of topics, which illustarate this basic idea
on the simple and tame examples without superfluous difficulties and details
as well as in the series [1] but from a bit more geometric point of view.

On the other hand this concrete article is placed at the crossing of two very
different ideologies of hidden symmetries (however, the ideological differences
may be rather subtle in practice). The used version of the first was
developed by G.Segal [6], M.Kontsevich, K.Gawedzki [7], M.Atiyah, G.Moore and
N.Seiberg [8], Yu.A.Neretin [9; and refs wherein], E.Witten [10] and others
[11]. This ideology, which was formulated most clearly in purely mathematical
fashion by G.Segal (pioneered the considered version), E.Witten, M.Atiyah and
to a certain extent by Yu.A.Neretin, may be characterized as a formal search of
hidden symmetries on the abstract level and is related to the direct problems of
representation theory. This ideology underlies the approaches of J.Mickellson
and D.P.Zhelobenko [12; and refs wherein] to the representation theory of
reductive Lie algebras and partially penetrates the researches of V.G.Drinfeld
and his numerous successors on Hopf and quasi-Hopf algebras. The second
manifested itself in the research activity of many specialists in mathematical
physics (e.g. in the Leningrad mathematical school's activity on the quantum
inverse scattering method [13; and refs wherein], or in the V.P.Maslov's group
investigations on nonlinear classical and quantum brackets [14,15; and refs
wherein]), the used below more formally mathematical version was elaborated by
the author [3]. The ideology may be characterized as a search of hidden
symmetries on the concrete level and is related to the inverse problems of
representation theory.

I hope that the presented crossing of two very different ideologies but
treating one subject will be rather interesting and crucial for the
understanding of various recently formed ``new looks'' in the representation
theory.

\head Topic Three: Projective modular functors and related objects
\endhead

\subhead 3.1. The projective semigroups $\Mantle(\PSLTWO)$ and
$\Voile(\PSLTWO)$: the mantle and veil (voile) of the projective group
$\PSLTWO$ [16:App.C.1]\endsubhead

The projective semigroup $\Mantle(\PSLTWO)$, the {\it mantle of the projective
group\/} $\PSLTWO$, can be realized in one of the following two ways. {\bf
Realization 1.} The elements of the semigroup $\Mantle(\PSLTWO)$ are
linear-fractional mappings $f:z\mapsto(az+b)(cz+d)^{-1}$ such that
$f(D_+)\subset D_+$, $f'\ne 0$. Multiplication of elements is a composition
of mappings. {\bf Realization 2.} The elements of the semigroup
$\Mantle(\PSLTWO)$ are domains $K$ in $\Bbb C$ that are homeomorphic to a
annulus and for which the components $\partial K_+$ and $\partial K_-$ of the
boundary $\partial K$ are circles in $\Bbb C$ with linear-fractional
parametrization. The inner circle $\partial K_+$ is parametrized in such a
way that as one passes around it the domain is on the right (such
parametrization is called ingoing), while for the outer circle $\partial K_-$
the domain is on the left (outgoing parametrization). Two domains are
said to be equivalent if there exists a linear-fractional automorphism of
$\bar\Bbb C$ that carries one of them into the other, with allowance for the
parametrization. Multiplication of elements is their glueing along the
parametrized boundaries.

The equivalence of the constructions is established as follows. If $f$ is a
linear-fractional mapping from $D_+$ to $D_+$, then as domain $K$ one may
consider the annulus $D_+\setminus\overset\circ\to D_+$, whose outer boundary has
the standard parametrization, while the parametrization of the inner boundary
is given by the mapping $f$. Conversely, if $K$ is an arbitrary domain that
satisfies the required conditions, then there exists a unique
linear-fractional automorphism of $\bar\Bbb C$ that maps $K$ to an annulus
whose outer boundary is the unit circle with the standard parametrization.
The required mapping $f$ is determined by the parametrization of the inner
component of the boundary.

The projective semigroup $\Mantle(\PSLTWO)$, the mantle of the projective
group $\PSLTWO$ is a partial complexification of it (cf.[17] and also [18]) and
has the complex dimension 3. An arbitrary finite-dimensional representation or
infinite-dimensional representation with highest weight of the projective
group $\PSLTWO$ can be extended by holomorphicity to a representation of
its mantle. This situation does not change if one consider projective
representations of both objects or linear representations of their universal
covers. Thus, the theory of representations of the projective group is,
in the words of Yu.A., Neretin, the theory of representations of a larger
semigroup that is ``invisible to the unaided eye'' -- its mantle, the
projective semigroup $\Mantle(\PSLTWO)$. This reformulation of the theory
of representations of $\PSLTWO$, which appears little more than tautologous
(though providing us with new interesting formulas or interpretations), has
great interest in that it stimulates one to look for further hidden structures
``invisible to the unaided eye'', for which the theory of representations is
richer than the theory of representations of the original projective group.
In particular, it would desirable if the theory of the irreducible
representations of this hidden structure automatically includes, in addition
to the theory of the irreducible representations of the projective group, the
corresponding Clebsch-Gordan coefficient calculus. The existence of such a
structure was indicated for a long time by the presence of a certain similarity
at the level of concrete formulas between the various objects of the theory
of irreducible representations and the Clebsch-Gordan coefficient calculus.
As an example, we can give the similarity of the matrix elements of
the irreducible representations and the corresponding Clebsch-Gordan
coefficients (see e.g.[19]).

As first step, we note that the construction of the semigroup
$\Mantle(\PSLTWO)$ can be generalized to arbitrary Riemann surfaces. The
semigroup $\Voile(\PSLTWO)$ that is then obtained is called the {\it veil\/}
({\it voile\/}) {\it of the projective group\/} $\PSLTWO$. The elements of
the projective semigroup $\Voile(\PSLTWO)$ are triples $(R,f_+,f_-)$,
where $R$ is a Riemann surface with the fixed projective structure, and
$f_+:D_+\hookrightarrow R$, $f_-:D_-\hookrightarrow R$ ($D_+=\{z: |z|\le 1\}$,
$D_-=\{z: |z|\ge 1\}$) are holomorphic imbeddings of the complex disks
$D_+$ and $D_-$ into $R$ (with allowance for the projective structure)
with noninterscting images. Multiplication of elements of the semigroup
is a sewing and so is defined in the similar way as multiplication in
Realization 2 of the projective semigroup $\Mantle(\PSLTWO)$.

The semigroup $\Voile(\PSLTWO)$ is a $\Bbb Z_+$-graded infinite-dimensional
semigroup. the grading is specified by the genus of the Riemann surface.
Formally, one is able to construct the noncommutative Grothendieck
group $\Gamma(\Voile(\PSLTWO))$ of the semigroup $\Voile(\PSLTWO)$.
The Grothendieck group $\Gamma(\Voile(\PSLTWO))$ is an infinite dimensional
group, however, its structure was not investigated.

The theory of representations of the projective semigroup $\Voile(\PSLTWO)$,
the veil of the projective group $\PSLTWO$, is much richer than the theory of
representations of its mantle. The number of representations depends on the
topology introduced on the semigroup. Among all topologies, the most
interesting are the following two: the ordinary topology on components of
fixed grading which discretely distinguishes the components, and a topology
that  takes into account possible continuous changes of genus. Problems of
the theory of representations of the semigroup $\Voile(\PSLTWO)$ as well
as its Grothendieck group $\Gamma(\Voile(\PSLTWO))$ are very important
but completely unexplored.

Some remarks should be added. The projective semigroup
$\Voile(\PSLTWO)$, the veil of the projective group $\PSLTWO$, can be
regarded as a ``fluctuating'' exponential of the Lie algebra $\sltwo$ that
describes the process of evolution with continuous creation and annihilation
of ``virtual particles''. From this point of view, it is very interesting
to study fractal (corresponding to nonperturbative effects) generalizations
of the veil $\Voile(\PSLTWO)$ of the projective group $\PSLTWO)$, the elements
of which can have infinite genus (cf.[20]). The corresponding completions
of the Grothendieck group $\Gamma(\Voile(\PSLTWO))$ are also of interest;
probably, some of such completions coincide with certain versions of one of
the infinite dimensional classical groups. Note that the nonperturbative
effects are essential for the quantum theory for a self-interacting string
field [21; and refs wherein].

\subhead 3.2. The manifold of projective trinions $\Trinion(\PSLTWO)$ and
its representations. Projective vertices [16:App.C.2]\endsubhead

Our next step will be to adapt the ideology of trinions [8] to the discussed
case.

A {\it projective trinion\/} is a quadruplet $(R,\partial R^1_+,\partial R^2_+,
\partial R_-)$, where $R$ is a Riemann surface of genus 0 ($R\subset\Bbb C$)
equipped with a projective structure and with a boundary whose components
$\partial R^1_+$, $\partial R^2_+$, $\partial R_-$ are homeomorphic to the
circle $S^1$ with ingoing linear-fractional parametrization defined on
$\partial R^1_+$ and $\partial R^2_+$ and outgoing linear-fractional
parametrization defined on $\partial R_-$ (the existence of such
parametrizations means that the real projective structures on the components
of the boundary, which are the restrictions of the complex projective structure
on the surface, are canonical).

The manifold of projective trinions $\Trinion(\PSLTWO)$ has the complex
dimension 6. The Lie group $G=\PSLTWO\times\PSLTWO\times\PSLTWO$ acts on the
set of trinions $\Trinion(\PSLTWO)$, and the corresponding action of the Lie
algebra $\frak g^{\Bbb C}=\sltwo+\sltwo+\sltwo$ is transitive. The stabilizer
of the trinion in the Lie algebra $\frak g^{\Bbb C}$ is the subalgebra of
holomorphic vector fields that admit extension to the complete trinion. Note
that the action of the Lie algebra $\frak g^{\Bbb C}$ on the manifold
$\Trinion(\PSLTWO)$ of projective trinions can be extended to the action of
the semigroup $\Mantle(\PSLTWO)\times\Mantle(\PSLTWO)\times\Mantle(\PSLTWO)$
with two copies of the projective semigroup $\Mantle(\PSLTWO)$, the mantle of
the projective group $\PSLTWO$, acting from the right, and one from the left.
At the same time, for all $\Mantle(\PSLTWO)$-integrable Verma modules $V_{h_1}$,
$V_{h_2}$, $V_{h_3}$ over the Lie algebra $\sltwo$ there exists no more than
one projective representation of the manifold of projective trinions
$\Trinion(\PSLTWO)$ (for the definition of a representation of a homogeneous
space based on the concept of a Mackey imprimitivity system, see [22,23]) in the
projective space $\BP(\Hom(V_{h_1}\otimes V_{h_2}; V_{h_3}))$ consistent with
the action of the projective semigroup $\Mantle(\PSLTWO)$, the mantle of the
projective group $\PSLTWO$, in these modules. Note that one may consider as
linear as projective representations of the $\Mantle(\PSLTWO)$, in the least
case one should consider the universal cover of the manifold
$\Trinion(\PSLTWO)$.

We now define the operation of {\it inserting a vertex\/} into an element of
the mantle $\Mantle(\PSLTWO)$. Let $g=(R,\partial R^1_+,\partial R^2_+,
\partial R_-)$ be a projective trinion, and $V_{h_1}$, $V_{h_2}$, $V_{h_3}$
be three integrable Verma modules. We denote by $A_g:V_{h_1}\otimes V_{h_2}
\mapsto V_{h_3}$ the operator (defined up to a factor) corresponding to the
projective trinion $g$ in the projective representation of the manifold
$\Trinion(\PSLTWO)$ of projective trinions in $\BP(\Hom(V_{h_1}\otimes V_{h_2};
V_{h_3}))$. Let $v$ be the highest vector in the Verma module $V_{h_1}$,
$k=(K,\partial K_+,\partial K_-)$ be an element of the projective semigroup
$\Mantle(\PSLTWO)$, the mantle of the projective group $\PSLTWO$, where
$K=R\sqcup_{\partial R^1_+}D$, $\partial K_+=\partial R^2_+$,
$\partial K_-=\partial R_-$ ($D$ is the disk bounded by $\partial R^1_+$), and
$u$ be an arbitrary point in $D$. Insertion of a vertex in the element $k$
means specification of an operator $A_k(u;v)$, defined up to a factor, in
$\Hom(V_{h_2};V_{h_3})$ as the limit of the family of operators $A_g$ as
$D\mapsto u$.

By means of the vertex insertion operation, we can define a projective
vertex. We consider an arbitrary element $k$ of the projective semigroup
$\Mantle(\PSLTWO)$, the mantle of the projective group $\PSLTWO$, with the
veretx $v$ of the weight $\mu$ that is inserted at the point $u$ and acts as
the operator $A_k(u;v)$ from the Verma module $V_h$ to the Verma module
$V_g$. A {\it projective vertex\/} $V_\mu(u;v)$ ($u\in S^1$) is the limit of
the family of operators $A_k(u;v)$ as $\partial K_+,\partial K_-\mapsto S^1$
(together with the parametrizations). The operator field $V_\mu(u;v)$ obtained
by means of this construction may be extended holomorphically with respect to
$u$. The result is just the $\sltwo$--primary field of weight $\mu$ from the
Verma module $V_h$ to the Verma module $V_g$. The definition and explicit
formulas for the $\sltwo$--primary fields in the Verma modules over the Lie
algebra $\sltwo$ may be found in [24]. Briefly, the $\sltwo$--primary fields
of weight (spin) $\mu$ may be characterized as generating functions for the
tensor operators, which transform as analytical $\mu$--differentials on a
circle $S^1$ (perhaps, with non-trivial monodromy).

\subhead 3.3. Polynions, their representations, and (derived) QPFT-operator
algebras [16:App.C.3;24]
\endsubhead

A {\it projective polynion\/} of degree $n$ is a data $(R,\partial R^1_+,
\partial R_+^2,\ldots\partial R_+^{n+1},\partial R_-)$, where $R$ is a Riemann
surface of genus 0 ($R\subset\bar\Bbb C$) that is equipped with a projective
structure and has a boundary whose components $\partial R^1_+$,
$\partial R^2_+$, ... $\partial R^{n+1}_+$, $\partial R_-$ are homeomorphic
to the circle $S^1$ with ingoing linear-fractional parametrization defined on
$\partial R^1_+$, $\partial R^2_+$, ... $\partial R^{n+1}_+$, and outgoing
linear-fractional parametrization on $\partial R_-$.

On the set of projective polynions $\Polynion(\PSLTWO)$ there are defined
the sewing operations
$$s:\Polynion(\PSLTWO)\times\Polynion(\PSLTWO)\mapsto\Polynion(\PSLTWO),$$
which are consistent with the grading
$$s:\Polynion_n(\PSLTWO)\times\Polynion_m(\PSLTWO)\mapsto
\Polynion_{n+m}(\PSLTWO).$$

The manifold $\Polynion_n(\PSLTWO)$ has the complex dimension $3(n+1)$. We
have
$$\aligned
\Polynion_0(\PSLTWO)&\simeq\Mantle(\PSLTWO),\\
\Polynion_1(\PSLTWO)&\simeq\Trinion(\PSLTWO),
\endaligned
$$
the polynions of degree greater than 1 being represented as compositions of
trinions. The Lie group $[\PSLTWO]^{n+2}$ acts on $\Polynion_n(\PSLTWO)$, and
the corresponding action of the Lie algebra $\frak g^{\Bbb C}=(n+2)\sltwo$ is
transitive. The stabilizer of a polynion in the Lie algebra $\frak g^{\Bbb C}
=(n+2)\sltwo$ is the subalgebra of holomorphic vector fields that admit
extension to the complete polynion. Note that the action of $\frak g^{\Bbb C}$
on $\Polynion_n(\PSLTWO)$ can be exponentiated to the action of
$[\Mantle(\PSLTWO)]^{n+2}$ with $n+1$ copies of the projective semigroup
$\Mantle(\PSLTWO)$, the mantle of the projective group $\PSLTWO$, acting from
the right, and one from the left.

The (projective) representation of polynions is a family of representations
$\pi_n$ of the homogeneous manifolds $\Polynion_n(\PSLTWO)$
$$\pi_n:\Polynion_n(\PSLTWO)\mapsto\BP(\Hom(V^{\otimes(n+1)};V))$$
such that the diagram
$$
\CD
\Polynion_n(\PSLTWO)\times\Polynion_m(\PSLTWO)@>>s>\Polynion_{n+m}(\PSLTWO)\\
@VV{\pi_n\times\pi_m}V                             @VV{\pi_{n+m}}V\\
\BP(\Hom(V^{\otimes(n+1)};V)\times\BP(\Hom(V^{\otimes(m+1)};V)@>>>
\BP(\Hom(V^{\otimes(n+m+1)};V))
\endCD
$$
where the lower arrow is the contraction operation, is commutative.

Let us now consider the infinitesimal counterparts of the representations
of projective polynions [25;16:\S 1] (see also [26] and refs wherein).
Some preliminary general concepts are necessary.

The {\it operator algebra of a quantum field theory\/} ({\it QFT-operator
algebra\/}) is a pair $(H,T^k_{ij}(\bold x))$ ($x$ belongs to $\Bbb R^n$ or
to $\Bbb C^n$), where $h$ is a linear space, and $t^k_{ij}(\bold x)$ is a
$H$-valued tensor field that satisfies $t^l_{im}(\bold x)t^m_{jk}(\bold y)=
t^m_{ij}(\bold x-\bold y)t^l_{mk}(\bold y)$. One may introduce the operators
$l_{\bold x}(e_i)$ ($e_i$ is an element of a basis in the space $H$):
$l_{\bold x}(e_i)e_j=t^k_{ij}(\bold x)e_k$. These operators satisfy the
identities $l_{\bold x}(e_i)l_{\bold y}(e_j)=t^k_{ij}(\bold x-\bold y)
l_{\bold y}(e_k)$ ({\it the operator product expansion\/}) and
$l_{\bold x}(e_i)l_{\bold y}(e_j)=l_{\bold y}l_{\bold x-\bold y}(e_i)e_j)$
({\it the duality\/}). Also one may define the multiplication operation
$m_{\bold x}$, which depends on the parameter $\bold x$, in the space
$H$: $m_{\bold x}(a,b)=l_{\bold x}(a)b$. For this operation {\it the identity
of smeared associativity\/} $m_{\bold x}(a,m_{\bold y}(b,c))=m_{\bold
y}(m_{\bold x-\bold y}(a,b),c)$ holds. The operators $l_{\bold x}$ are the
operators of multplication from the left in the obtained algebra, and
$t^k_{ij}(\bold x)$ are the structural functions. Such definition of the
QFT-operator algebras is an axiomatization of the well-known operator product
expansions in quantum field theory.

In a QFT-operator algebra with unit there is defined the operator
$\bold L$: $\bold La=\frac{d}{d\bold x}\left.(l_{\bold x}(a)\boldkey 1)
\right|_{\bold x=0}$. This operator generates the infinitesimal translations
$[\bold L,l_{\bold x}(a)]=\frac{d}{d\bold x}l_{\bold x}(a)$. In what follows,
we shall assume that the variable $\bold x$ ranges over the complex plane and
that the tensor field $t_{ij}^k(\bold x)$ is analytic.

Let us now describe an object, which is an infinitesimal counterpart of the
representations of polynions. A QFT-operator algebra $(H,t_{ij}^k(u))$ is
called a {\it (derived) QPFT-operator algebra\/} (and, by some authors,
quasi-vertex algebra) iff
\roster
\item"(1)" the linear space $H$ can be decomposed into a direct sum or a
direct integral of Verma modules $V_{\alpha}$ over the Lie algebra
$\sltwo$ with the highest vectors $v_{\alpha}$ and the highest weight
$h_{\alpha}$;
\item"(2)" the operator fields $l_u(v_{\alpha})$ are $\sltwo$--primary
fields of weight $h_{\alpha}$, in other words, on commutation with the
generators of the Lie algebra $\sltwo$ they transform in accordance with a
tensor law as $h_{\alpha}$--differentials:
$$[L_k,l_u(v_{\alpha})]=(-u)^k(u\partial_u+(k+1)h_{\alpha})l_u(v_{\alpha}),
\qquad k=-1,0,1;$$
\item"(3)" the following derivative rule of generation of descendents holds:
$$[L_{-1},l_u(f)]=\tfrac{d}{du}l_u(f)=l_u(L_{-1}f).$$
\endroster

Let us formulate now the main result establishing a connection between the
representations of polynions and (derived) QPFT-operator algebras. {\sl
If the representation space $V$ of the polynions is decomposed into a direct
sum or a direct integral of the Verma modules over the Lie algebra $\sltwo$,
then the transition to projective vertices in the representation of polynions
defines the structure of a (derived) QPFT-operator algebra in the
representation space}. In general, the converse is not true; not every (derived)
QPFT-operator algebra can be integrated to a representation of polynions (in
the same way that not every representation of a Lie algebra can be integrated
to a representation of the corresponding Lie group). One may consider the
relations between (derived) QPFT-operator algebras and representations
of polynions as analogous to ones between Lie algebras and Lie groups
[16:App.C].

\remark{Remark: Local projective field algebras}
The method of smoothing (that means the transition from the operators
representing elements of a group to their integrals representing elements
of the group algebra) is very effective in the representation theory (see
e.g.[22,27]).
The analogous procedure may be applied to the QFT-operator algebras
[28] (one should use the smoothing by a ``vertex position'' $\bold x$).
In the case of the QPFT-operator algebras the result may be defined
axiomatically; the obtained object is called a {\it local projective field
algebra\/} [26], it may be considered as a ``vertex analogue'' of the group
algebra. Note that one may construct an analogue of the group algebra
(``polynionic algebra'') by an integration (``smoothing'') of the operators
representing polynions.
\endremark

\subhead 3.4. The projective category $\Train(\PSLTWO)$, the train of the
projective group $\PSLTWO$, and the projective modular functor [16:App.C.4]
\endsubhead

The projective category $\Train(\PSLTWO)$, the train of the projective group
$\PSLTWO$, is a category whose objects $\frak A$, $\frak B$, $\frak C$, ...
are finite ordered sets, morphisms in $\Mor(\frak A,\frak B)$ are the data
$(R,\partial R^1_+,\ldots\partial R^{n_+}_+,\partial R^1_-,\ldots\partial
R^{n_-}_-, n_+=\#\frak A, n_-=\#\frak B)$, where $R$ is a Riemann surface
equipped with projective structure and possessing a boundary whose components
$\partial R^i_+$, $\partial R^i_-$ are homeomorphic to the circle $S^1$ with
ingoing linear-fractional parametrization on $\partial R^1_+,\ldots\partial
R^{n_+}_+$, $n_+=\#\frak A$, and outgoing linear-fractional parametrization
on $\partial R^1_-,\ldots \partial R^{n_-}_-$, $n_-=\#\frak B$. Composition
of morphisms is the sewing operation $s$.

On the set $\Mor(\frak A,\frak B)$ there acts the Lie group
$[\PSLTWO]^{\#\frak A+\#\frak B}$, but in contrast to polynions the
corresponding action of the Lie algebra $\frak g^{\Bbb C}=
(\#\frak A+\#\frak B)\sltwo$ is not transitive (this being due to the
presence of moduli of Riemann surfeces of nonvanishing genus). The stabilizer
of a  morphism in the Lie algebra $\frak g^{\Bbb C}$ is the subalgebra of
holomorphic vector fields that admit an extension to the geometrical image
of the morphism. The action $\frak g^{\Bbb C}$ on the set of morphisms can
be exponentiated to the action of the semigroup
$[\Mantle(\PSLTWO)]^{\#\frak A+\#\frak B}$, with $\#\frak A$ copies of the
projective semigroup $\Mantle(\PSLTWO)$, the mantle of the projective group,
acting from the right, and $\#\frak b$ from the left. Although the action of
the Lie algebra $\frak g^{\Bbb C}$ on the set of morphisms is not transitive,
it is possible to define the concept of a representation of the family of
morphisms as a continuous family of representations of the orbits of this
action.

A (projective) representation of the projective category $\Train(\PSLTWO)$,
the train of the projective group $\PSLTWO$ ({\it projective modular
functor\/}), is a set of representations $\pi_{\frak A\frak B}:\Mor(\frak
A,\frak B)\mapsto\BP(\Hom(V^{\otimes\#\frak A},V^{\otimes\#\frak B}))$ such
that
\roster
\item"--" the diagram
$$
\CD
\Mor(\frak A,\frak B)\times\Mor(\frak B,\frak C)@>>s>\Mor(\frak A,\frak C)\\
@VV\pi_{\frak A\frak B}\times\pi_{\frak B\frak C}V @VV\pi_{\frak A\frak C}V\\
\BP(\Hom(V^{\otimes\#\frak A},V^{\otimes\#\frak B})\times
\BP(\Hom(V^{\otimes\#\frak B},V^{\otimes\#\frak C}) @>>>
\BP(\Hom(V^{\otimes\#\frak A},V^{\otimes\#\frak C})),
\endCD
$$
where the lower arrow is the contraction operation, is commutative;
\item"--" if $\frak A=\frak A_1\sqcup\frak A_2$, $\frak B=\frak B_1\sqcup
\frak B_2$, $R=R_1\sqcup R_2$, $R_i\in\Mor(\frak A_i,\frak B_i)$, then
$\pi_{\frak A\frak B}(R)=\pi_{\frak A_1\frak B_1}(R_1)\times
\pi_{\frak A_2\frak B_2}(R_2)$.
\endroster

Every projective modular functor corresponds to some representation of
poly\-nions, since polinions are a special case of morphisms in the
projective category $\Train(\PSLTWO)$, the train of the projective group
$\PSLTWO$. In general, the converse is not true -- not every representation
of polynions can be extended to a projective modular functor. Indeed, the
projective semigroup $\Voile(\PSLTWO)$, the veil of the projective group
$\PSLTWO$, is identified with the semigroup of all endomorphisms of object of
cardinality 1 in the projective category $\Train(\PSLTWO)$, and this is the
``topological'' obstruction to an extension of representations of polynions to
projective modular functors. Thus, among the structures of the theory of
represesentations in {\it the quantum projective field theory\/} (derived
QPFT-operator algebras, representations of projective polynions, projective
modular functors) the last [the representations of the projective category
$\Train(\PSLTWO)$, the train of the projective group $\PSLTWO$] form the
smallest class (under the condition that the representation spaces are sums of
Verma modules over the Lie algebra $\sltwo$), and the first [the derivative
QPFT-operator algebras or, equivalently, the QPFT-operator algebras] form
the largest class.

\remark{Remark: ``Pseudogroup'' variations on the theme of ``Train''}
Let $\Cal N$ be any finite set. Let us consider some lattice $\Cal P$ of
subsets of $\Cal N$ ($\Cal P$ is closed under intersections and unions). One
may construct a generalization $\Train(\PSLTWO,\Cal P)$ of the projective
category $\Train(\PSLTWO)$, which objects belong to $\Cal P$. Morphisms are
defined in the same manner as for $\Train(\PSLTWO)$. The category
$\Train(\PSLTWO,\Cal P)$ is supplied with a natural structure of a topologized
category.

Certainly, one may consider the manifold $\widetilde{\Cal N}=\Cal N\times S^1$
with the projective structure instead of $\Cal N$. Some elements of the
category $\Train(\PSLTWO,\Cal P)$ will form a pseudogroup of projective
transformations of $\widetilde{\Cal N}$. Though other elements of
$\Train(\PSLTWO,\Cal P)$ do not constitute a pseudogroup, its categorical
properties are analogous to ones of ordinary pseudogroups and may be
straightforwardly axiomatized. Some examples of analogous (``pseudogroup'')
categories were considered by Yu.A.Neretin [9]. The representations of such
``pseudogroup'' categories may be naturally defined in the same universal
manner. Note that the representations of smooth pseudogroups of transformations
appear in the framework of the asymptotic quantization [15:Ch.4]. The
algorythm of asymptotic quantization uses the representations ``mod $\hbar$''
($\hbar$ is a parameter). The infinitesimal counterparts of such
representations, the asymptotic representations of the pseudoalgebras of
vector fields ``mod $\hbar$'', are a partial case of the general
$\frak A$--projective representations of the topic 10 of series [2].
\endremark

\remark\nofrills{Remark:}\ In all constructions above the real projective
structures on the boundaries of the suitable Riemann surfaces were canonical
so that the components of the boundaries admitted linear-fractional
parametrizations. One may generalize the situation omitting this condition.
\endremark

\subhead 3.5. Conclusions\endsubhead

Thus, a general scheme of the reconstruction of hidden objects related to
the Lie groups is briefly exposed above on the simplest example of
the projective group $\PSLTWO$. However, it may be evidently generalized to
other finite-dimensional semisimple Lie groups with changes in minor details.
In the next topic we shall discuss how it is adapted to the
infinite-dimensional group $\Diff_+(S^1)$ of all orientation preserving
diffeomorphisms of a circle $S^1$ (or its central extension $\Vir$, the
Virasoro-Bott group). The exposition of the general scheme above lacks
a lot of interesting details (such as QPFT-operator crossing-algebras,
which are the infinitesimal counterparts of projective modular functors,
and many related structures or the projective Krichever-Novikov functors,
cf.[21]) for the simplicity and clarity. The more detalized exposition of
the scheme should be found in [21]. The relation of the reconstructed
algebraic objects to {\sl the quantum group foundations of the
self-interacting string field theory\/} is also described in the article [21],
which contains a necessary bibliography. Here we mark
only that the (renormalized) version of the projective Krichever-Novikov
functor is based on the concept of a sewing of noncommutative coverings of
the Riemann surfaces [29], which are realized by means of sheaves of the
local projective field algebras. The original version of the Krichever-Novikov
construction is adapted to the conformal case [30] and has a deal with
operator product expansions on the Riemann surfaces instead of local field
algebras.

Note that our scheme is ideologically the same as one of G.Segal, E.Witten and
M.Atiyah but differs from the essentially less known general scheme of
Yu.A.Neretin, who systematically avoids, neglects or reduces any topological
effects in the reconstruction of the hidden objects (at least, such effects
appear only episodically in his articles and in my opinion their appearing is
motivated presumably by influences of the other authors such as G.Segal,
M.Kontsevich, M.Atiyah, I.M.Krichver and S.P.Novikov, G.Moore or E.Witten),
and in such a way describes some very interesting purely algebraic or analytic
phenomena [9; and refs wherein]. Note that the topological aspects of an
analysis of hidden symmetries in the representation theory connect the least
with {\sl the bordism theory\/} (indeed, the veil $\Voile(\PSLTWO)$ of the
projective group $\PSLTWO$ is just the semigroup of all bordisms of a circle
supplied with the canonical projective structure).

Here it is convenient to formulate some open problems related to the
material of this topic.

\remark\nofrills{Problems:}
\roster
\item"--" To generalize the scheme to the simplest nonlinear objects such as
the Racah-Wigner algebra or the higher $\mho$--algebras for $\sltwo$ [2;3:\S1]
and the Sklyanin algebra [31;3:\S3.1]. It is especially interesting to define
the ``fluctuating exponents'' for the nonlinear objects similar to the
Racah-Wigner algebras. Such construction should be considered in a general
context of the nonlinear geometric algebra [32] and quantization of nonlinear
Poisson brackets [15].
\item"--" To generalize the scheme to the isotopic pairs [33;3:\S2.3],
especially to R-matrix ones.
\item"--" To generalize the scheme to the ``quantum $\sltwo$'' [34] and,
perhaps only partially, to other ``nonlinear $\sltwo$'' [35].
\endroster
\endremark

\head Topic Four: Conformal modular functors and related objects
\endhead

The objects of this topic are infinite-dimensional analogs of ones
discussed above, and so they constitute a subject of the representation
theory of the infinite dimensional Lie groups, Lie algebras and related
structures. Note that the procedures, which were almost trivial and
tautological for the finite-dimensional case, become very profound in its
infinite-dimensional counterpart. For instance, the construction of mantles
of the infinite-dimensional groups is a part of a general ideology of
G.I.Olshanski{\v\i} [36] of the semigroup approach to the representation
theory of such groups.

\subhead 4.1. The Lie algebra $\Vect(S^1)$ of vector fields on a circle,
the group $\Diff_+(S^1)$ of diffeomorphisms or a circle, the Virasoro
algebra $\vir$, the Virasoro-Bott group $\Vir$ and the Neretin semigroup
$\Ner$
\endsubhead

Let $\Diff(S^1)$ denote the group of diffeomorphisms of the unit circle
$S^1$. The group manifold $\Diff(S^1)$ splits into two connected components,
the subgroup $\Diff_+(S^1)$ and the coset $\Diff_-(S^1)$. The diffeomorphisms
in $\Diff_+(S^1)$ preserve the orientation on the circle $S^1$ and those in
$\Diff_-(S^1)$ reverse it.

The Lie algebra of $\Diff_+(S^1)$ can be identified with the linear space
$\Vect(S^1)$ of smooth vector fields on the circle equipped with the
commutator
$$[v(t)d/dt,u(t)d/dt]=(v(t)u'(t)-v'(t)u(t))d/dt.$$
In the basis $s_n=\sin(nt)d/dt$, $c_n=\cos(nt)d/dt$, $h=d/dt$ the
commutation relations have the form
$$\aligned
[s_n,s_m]&=0.5((m-n)s_{m+n}+\sgn(n-m)(n+m)s_{|n-m|}),\\
[c_n,c_m]&=0.5((n-m)s_{n+m}+\sgn(n-m)(n+m)s_{|n-m|}),\\
[s_n,c_m]&=0.5((m-n)c_{n+m}-(n+m)c_{|n-m|})-n\delta_{nm}h,\\
[h,s_n]&=nc_n,\quad [h,c_n]=ns_n.
\endaligned
$$
The generators $h$, $s_n$, $c_n$ form a Lie algebra isomorphic to
$\sltwor$ for each $n$.

The complexification of the Lie algebra $\Vect(S^1)$ will be denoted by
$\CVect(S^1)$. It is convenient to choose the basis $e_k=ie^{ikt}d/dt$ in
$\CVect(S^1)$. The commutation relations of the algebra $\CVect(S^1)$ have
the following form
$$[e_j,e_k]=(j-k)e_{j+k}$$
in the basis $e_k$. The generators $e_n$, $e_{-n}$, $e_0$ form a Lie algebra
isomorphic to $\sltwo$ for each $n$.

In 1968 I.M.Gelfand and D.B.Fuchs discovered [37] that $\Vect(S^1)$ possesses
a nontrivial central extension. The corresponding 2-cocycle is
$c(u,v)=\int v'(t)du'(t)$ or, equivalently, $c(u,v)=\left|\matrix
v'(t_0) & u'(t_0) \\ v''(t_0) & u''(t_0) \endmatrix\right|$ (see [38]).
This central extension was independently discovered by M.Virasoro [39] and
named after him. Let us denote the Virasoro algebra by $\vir$. Its
complexification, which is also called the {\it Virasoro algebra}, will be
denoted $\Cvir$. As a vector space $\vir$ is generated by the vectors $e_k$
and the central element $c$. The commutation relations have the form
$$[e_j,e_k]=(j-k)e_{j+k}+\delta(j+k)\tfrac{j^3-j}{12}c.$$
The imbeddings of the Lie algebras $\sltwor$ and $\sltwo$ into $\Vect(S^1)$
and $\CVect(S^1)$ may be lifted to the imbeddings of these Lie algebras into
$\vir$ and $\Cvir$.

The infinite-dimensional group $\Vir$ corresponding to the algebra $\vir$ is
a central extension of the group $\Diff(S^1)$. The corresponding 2-cocycle
was calculated by R.Bott [40] so the group $\Vir$ is called the {\it
Virasoro-Bott group}. The imbeddings of the Lie algebra $\sltwor$ into
$\Vect(S^1)$ or $\vir$ are exponentiated to the imbeddings of the
$n$-coverings of the projective group $\PSLTWO$ into
$\Diff_+(S^1)$ and $\Vir$ ($n$ labels the imbedding).

In 1968 in the cited paper I.M.Gelfand and D.B.Fuchs computed completely
the cohomology of $\Vect(S^1)$ and, thus, discovered a non-trivial 3-cocycle
of this Lie algebra. It may be written as $c(u,v,w)=\det(B(t_0))$,
where
$$B(t)=\left[\matrix u'(t) & v'(t) & w'(t) \\ u''(t) & v''(t) &
w''(t) \\ u'''(t) & v'''(t) & w'''(t) \endmatrix\right],$$
or as $c(u,v,w)=\int\det(B(t))dt$ [38]. This cocycle defines [41:\S1.1] a {\it
nonassociative deformation\/} of the Virasoro-Bott group, which was called
the {\it Gelfand-Fuchs loop} and denoted by $\Gf$.

The are no groups corresponding to the Lie algebras $\CVect(S^1)$ or $\Cvir$,
but one can consider the following construction due to Yu.A.Neretin,
M.L.Kontsevich and G.Segal. Let us denote by $\LDiff^{\Bbb C}_+(S^1)$ the set
of all analytic mappings $g:S^1\mapsto\Bbb C\setminus\{0\}$ such that $g(S^1)$
is a Jordan curve surrounding zero, the orientations of $S^1$ and $g(S^1)$ are
the same, and $g'(e^{i\theta})$ is everywhere different from zero.
$\LDiff^{\Bbb C}_+(S^1)$ is a local group. Let $\LNer\subset\LDiff^{\Bbb
C}_+(S^1)$ be the local subsemigroup of mappings $g$ such that
$|g(e^{i\theta})|<1$. As it was shown by Yu.A.Neretin in 1987, the structure of
local semigroup on $\LNer$ extends to the structure of global semigroup $\Ner$.
There exist at least two constructions of the semigroup $\Ner$.

{\it The first construction} (Yu.A.Neretin). An element of $\Ner$ is a formal
product $p\cdot A(t)\cdot q$ (*), where $p,q\in\Diff_+(S^1)$, $p(1)=1$, $t>0$,
$A(t):\Bbb C\mapsto\Bbb C$, $A(t)z=e^{-t}z$. To define multiplication in $\Ner$,
one must describe the rule used to transform the formal product
$A(s)\cdot p\cdot A(t)$ to the form (*).

{\bf A.} Let $t$ be so small that the diffeomorphism $p$ extends
holomorphically to the annulus $e^{-t}\le|z|\le1$. Then the product $g=A(s)pA(t)$
is well defined. Let $K$ be the domain bounded by $S^1$ and $g(S^1)$. Let $Q$
be the canonical conformal mapping of $K$ onto the annulus $e^{-t'}\le z\le 1$,
normalized by the condition $Q(1)$. Then $g=p'\cdot A(t)\cdot q'$, where
$p'=\left.Q^{-1}\right|_{S^1}$ and $q'$ is determined by the identity
$A(s)\cdot p\cdot A(t)=p'\cdot A(t')\cdot q'$.

{\bf B.} For an arbitrary $t$ there exists a suitable $n$ such that the
product
$$A(s)\cdot p\cdot A(t)=(\ldots(A(s)\cdot p\cdot A(t/n))A(t/n)\ldots)A(t/n)
\tag{**}$$
can be calculated. It can be shown that the product does not depend on the
choice of the representation (**) and is associative.

{\it The second construction\/} (M.L.Kontsevich and G.Segal). An element $g$
of the semigroup $\Ner$ is a triple $(K,p,q)$, where $K$ is a Riemann surface
with boundary $\partial K$ such that $K$ is biholomorphically equivalent to
the annulus and $p,q:S^1\mapsto\partial K$ are fixed parametrizations of
the components of $\partial K$. The mapping $p$ realizes the ingoing
parametrization whereas $q$ realizes the outgoing parametrization. Two
elements $g_i=(K_i,p_i,q_i)$, $i=1,2$, are equivalent if there exists a
conformal mapping $R:K_1\mapsto K_2$ such that $p_2=Rp_1$, and $q_2=Rq_1$.
The product of two elements $g_1$ and $g_2$ is the element $g_3=(K_3,p_3,q_3)$,
where
$$K_3=K_1\bigsqcup_{q_1(e^{it}=p_2(e^{it})}K_2,$$
$p_3=p_1$, and $q_3=q_2$.

The Neretin semigroup $\Ner$ is called the mantle of the group $\Diff_+(S^1)$
of diffeomorphisms of a circle and is denoted by $\Mantle(\Diff_+(S^1))$
The Neretin semigroup $\Ner$ possesses a central extension, which was
discovered by Yu.A.Neretin in 1989. This extension is compatible with the
Virasoro-Bott extension $\Vir$ of the group $\Diff_+(S^1)$ and is a mantle
$\Mantle(Vir)$ of the group $\Vir$. The imbeddings of the $n$-coverings of
$\PSLTWO$ into the groups $\Diff_+(S^1)$ and $\Vir$ may be extended to the
holomorphic imbeddings of the mantles of such coverings (which are coverings
of $\Mantle(\PSLTWO)$) into the mantles of the infinite-dimensional groups
$\Diff_+(S^1)$ and $\Vir$. The first imbedding ($n=1$) may be used for
the construction of a certain enlargement of the semigroup
$\Mantle(\Diff_+(S^1))$, namely,
$\Mantle(\Diff_+(S^1))\times_{\Mantle(\PSLTWO)}\PSLTWOC$.

Note that in the infinite-dimensional case each projective representation
of the group $\Diff_+(S^1)$ is linearized over the universal covering of
its central extension $\Vir$ whereas in the finite-dimensional case
each projective representation of the group $\PSLTWO$ is linearized over
its own universal covering. The same situation is for the mantles.

The further interesting information on the objects of this paragraph and
their relation to the infinite-dimensional geometry should be found in the
papers [42,43].

\subhead 4.2. The semigroup $\Voile(\Diff_+(S^1))$ (the veil of the group
$\Diff_+(S^1)$), the manifolds $\Trinion(\Diff_+(S^1))$ of (conformal) trinions
and $\Polynion(\Diff_+(S^1))$ of (conformal) polynions. The representations
of polynions and QCFT-operator algebras [21]
\endsubhead

The construction of the semigroup $\Mantle(\Diff_+(S^1))$ can be generalized
to arbitrary Riemann surfaces. The semigroup $\Voile(\Diff_+(S^1))$ that is
then obtained is called the {\it veil\/} ({\it voile\/}) {\it of the group\/}
$\Diff_+(S^1))$. The elements of the semigroup $\Voile(\Diff_+(S^1))$ are
triples $(R,f_+,f_-)$, where $R$ is a Riemann surface, and
$f_+:D_+\hookrightarrow R$, $f_-:D_-\hookrightarrow R$ ($D_+=\{z: |z|\le 1\}$,
$D_-=\{z: |z|\ge 1\}$) are holomorphic imbeddings of the complex disks
$D_+$ and $D_-$ into $R$ with nonintersecting images. Multiplication of
elements of the semigroup is a sewing. The semigroup $\Voile(\Diff_+(S^1))$
is a $\Bbb Z_+$-graded infinite-dimensional semigroup. The grading is
specified by the genus of the Riemann surface. The semigroup
$\Voile(\Diff_+(S^1))$, the veil of the group $\Diff_+(S^1))$, can be
regarded as a ``fluctuating'' exponential of the Lie algebra $\Vect(S^1)$ that
describes the process of evolution with continuous creation and annihilation
of ``virtual particles''. From this point of view, it is very interesting
to study fractal (corresponding to nonperturbative effects) generalizations
of the veil $\Voile(\Diff_+(S^1))$ of the projective group $\Diff_+(S^1))$,
the element of which can have infinite genus; note once more that the
nonperturbative effects in this situation are essential for the quantum theory
for a self-interacting string field [21; and refs wherein].

Let us briefly repeat the construction of trinions [8].

A {\it (conformal) trinion\/} is a quadruplet $(R,\partial R^1_+,\partial R^2_+,
\partial R_-)$, where $R$ is a Riemann surface of genus 0 ($R\subset\Bbb C$)
with a boundary whose components $\partial R^1_+$, $\partial R^2_+$,
$\partial R_-$ are homeomorphic to the circle $S^1$ with ingoing
parametrization defined on $\partial R^1_+$ and $\partial R^2_+$ and outgoing
parametrization defined on $\partial R_-$.

The situation differs from one described in Topic 3 by the lack of projective
structure and the related conditions on the ingoing and outgoing
parametrizations.

The manifold of conformal trinions $\Trinion(\Diff_+(S^1))$ is an
infinite-dimensional complex manifold. The Lie group
$G=\Diff_+(S^1)\times\Diff_+(S^1)\times\Diff_+(S^1)$ acts on the set of
trinions $\Trinion(\Diff_+(S^1))$, and the corresponding action of the Lie
algebra $\frak g^{\Bbb C}=\CVect(S^1)+\CVect(S^1)+\CVect(S^1)$ is transitive.
The stabilizer of the trinion in the Lie algebra $\frak g^{\Bbb C}$ is the
subalgebra of all holomorphic vector fields that admit extension to the
complete trinion. Note that the action of the Lie algebra $\frak g^{\Bbb C}$
on the manifold $\Trinion(\Diff_+(S^1))$ of conformal trinions can be extended
to the action of the semigroup $\Mantle(\Diff_+(S^1))\times
\Mantle(\Diff_+(S^1))\times\Mantle(\Diff_+(S^1))$ with two copies of the
semigroup $\Mantle(\Diff_+(S^1))$, the mantle of the group $\Diff_+(S^1)$,
acting from the right, and one from the left. At the same time, for all Verma
modules $V_{h_1,c}$, $V_{h_2,c}$, $V_{h_3,c}$ over the Lie algebra $\Cvir$
[44] (here $h_i$ are the extremal weights and $c$ is the central charge), which
are integrable to the projective representations of the semigroup
$\Mantle(\Diff_+(S^1))$ there exists no more than one projective
representation of the manifold of conformal trinions $\Trinion(\Diff_+(S^1))$
in the projective space $\BP(\Hom(V_{h_1,c}\otimes V_{h_2,c}; V_{h_3,c}))$
consistent with the action of the semigroup $\Mantle(\Diff_+(S^1))$ in these
modules. One can also consider the universal cover of the manifold
$\Trinion(\Diff_+(S^1))$.

The operation of the vertex insertion and the (conformal) vertices themselves
can be defined in a way similar to one specified in the Topic 3.

A {\it (conformal) polynion\/} of degree $n$ is a data $(R,\partial R^1_+,
\partial R_+^2,\ldots\partial R_+^{n+1},\partial R_-)$, where $R$ is a Riemann
surface of genus 0 ($R\subset\bar\Bbb C$) that
has a boundary whose components $\partial R^1_+$,
$\partial R^2_+$, ... $\partial R^{n+1}_+$, $\partial R_-$ are homeomorphic
to the circle $S^1$ with ingoing parametrization defined on
$\partial R^1_+$, $\partial R^2_+$, ... $\partial R^{n+1}_+$, and outgoing
parametrization on $\partial R_-$.

The situation differs from one described in Topic 3 by the lack of projective
structure and the related conditions on the ingoing and outgoing
parametrizations.

On the set of conformal polynions $\Polynion(\Diff_+(S^1))$ there are defined
the sewing operations
$$s:\Polynion(\Diff_+(S^1))\times\Polynion(\Diff_+(S^1))\mapsto
\Polynion(\Diff_+(S^1)),$$
which are consistent with the grading
$$s:\Polynion_n(\Diff_+(S^1))\times\Polynion_m(\Diff_+(S^1))\mapsto
\Polynion_{n+m}(\Diff_+(S^1)).$$

The manifold $\Polynion_n(\Diff_+(S^1))$ is an infinite-dimensional complex
manifold. We have
$$\aligned
\Polynion_0(\Diff_+(S^1))&\simeq\Mantle(\Diff_+(S^1)),\\
\Polynion_1(\Diff_+(S^1))&\simeq\Trinion(\Diff_+(S^1)),
\endaligned
$$
the polynions of degree greater than 1 being represented as compositions of
trinions. The Lie group $[\Diff_+(S^1)]^{n+2}$ acts on
$\Polynion_n(\Diff_+(S^1))$, and the corresponding action of the Lie algebra
$\frak g^{\Bbb C}=(n+2)\CVect(S^1)$ is transitive. The stabilizer of a
polynion in the Lie algebra $\frak g^{\Bbb C}=(n+2)\CVect(S^1)$ is the
subalgebra of all holomorphic vector fields that admit an extension to the
complete polynion. Note that the action of $\frak g^{\Bbb C}$ on
$\Polynion_n(\Diff_+(S^1))$ can be exponentiated to the action of
$[\Mantle(\Diff_+(S^1))]^{n+2}$ with $n+1$ copies of the projective semigroup
$\Mantle(\Diff_+(S^1))$, the mantle of the group $\Diff_+(S^1)$, acting from
the right, and one from the left.

The (projective) representation of polynions is a family of representations
$\pi_n$ of the homogeneous manifolds $\Polynion_n(\Diff_+(S^1))$
$$\pi_n:\Polynion_n(\Diff_+(S^1))\mapsto\BP(\Hom(V^{\otimes(n+1)};V))$$
such that the diagram
$$
\CD
\Polynion_n(\Diff_+(S^1))\times\Polynion_m(\Diff_+(S^1))@>>s>
\Polynion_{n+m}(\Diff_+(S^1))\\
@VV{\pi_n\times\pi_m}V                             @VV{\pi_{n+m}}V\\
\BP(\Hom(V^{\otimes(n+1)};V)\times\BP(\Hom(V^{\otimes(m+1)};V)@>>>
\BP(\Hom(V^{\otimes(n+m+1)};V))
\endCD
$$
where the lower arrow is the contraction operation, is commutative.

The extremal vector $T$ of the weight 2 in the QPFT-operator algebra
is called a {\it conformal stress-energy tensor}, if $T(u):=l_u(T)=\sum L_k
(-u)^{k-2}$, and operators $L_k$ generate the Virasoro algebra:
$[L_i,L_j]=(i-j)L_{i+j}+\frac{i^3-i}{12} c\cdot I$.
In view of results of [45] the QPFT-operator algebras with the conformal
stress-energy tensor are just the {\it operator algebras of the quantum
conformal field theory\/} ({\it QCFT-operator algebras\/}) in sense of
the papers [46,29].

A connection between the representations of polynions and QCFT-operator
algebras was established in the article [21]. {\sl If the representation
space $V$ of the polynions is decomposed into a direct sum or a direct
ntegral of the Verma modules over the Lie algebra $\Cvir$, then the
transition to (conformal) vertices in the representation of polynions
defines the structure of a QCFT-operator algebra in the representation
space}. In general, the converse is not true; not every QCFT-operator algebra
can be integrated to a representation of polynions (in the same way that not
every representation of a Lie algebra can be integrated to the representation
of the corresponding Lie group). One may consider the relations between
QCFT-operator algebras and representations of polynions as analogous to
ones between Lie algebras and Lie groups.

\remark\nofrills{Remarks:}
\roster
\item"--" There is a subtle distinction between the derived
QPFT-operator algebras, which were defined above, and the QPFT-operator
algebras, which were considered in [46,29]. The fact that such distinction
is subtle indeed and may be neglected was explicated in the article [25].
\item"--" Some subclasses of the QCFT-operator algebras were considered in
[47; and refs wherein]  under the title of ``vertex operator algebras'' and
``vertex algebras'' (see [21]).
\item"--" The local projective field algebras for QCFT-operator algebras
(the local conformal field algebras) were defined and explored in
articles [48]. Their geometric interpretation (as structural rings of
the noncommutative coverings of Riemann surfaces) was given in [29].
The relations between the local conformal and projective field algebras
[the ``projective reduction''] were analized in [49] (see also [29]).
\item"--" One may consider the projective and conformal {\it antitrinions\/}
and {\it an\-ti\-po\-ly\-ni\-ons}, which are received from trinions and
polynions by the changing of the ingoing parametrizations to the outgoing ones
and vice versa, as well as their representations in the way analogous to one
described above (see [21]).
\endroster
\endremark

\subhead 4.3. The category $\Train(\Diff_+(S^1))$, the train of the
group $\Diff_+(S^1)$, and the (conformal) modular functor [21]
\endsubhead

The category $\Train(\Diff_+(S^1))$, the train of the group $\Diff_+(S^1)$,
is a category whose objects $\frak A$, $\frak B$, $\frak C$, ...
are finite ordered sets, morphisms in $\Mor(\frak A,\frak B)$ are the data
$(R,\partial R^1_+,\ldots\partial R^{n_+}_+,\partial R^1_-,\ldots\partial
R^{n_-}_-, n_+=\#\frak A, n_-=\#\frak B)$, where $R$ is a Riemann surface
possessing a boundary whose components
$\partial R^i_+$, $\partial R^i_-$ are homeomorphic to the circle $S^1$ with
ingoing parametrization on $\partial R^1_+,\ldots\partial
R^{n_+}_+$, $n_+=\#\frak A$, and outgoing parametrization
on $\partial R^1_-,\ldots \partial R^{n_-}_-$, $n_-=\#\frak B$. Composition
of morphisms is the sewing operation $s$.

The situation differs from one described in Topic 3 by the lack of projective
structure and the related conditions on the ingoing and outgoing
parametrizations.

On the set $\Mor(\frak A,\frak B)$ there acts the Lie group
$[\Diff_+(S^1)]^{\#\frak A+\#\frak B}$, but in contrast to polynions the
corresponding action of the Lie algebra $\frak g^{\Bbb C}=
(\#\frak A+\#\frak B)\Cvir$ is not transitive (this being due to the
presence of moduli of Riemann surfeces of nonvanishing genus). The stabilizer
of a  morphism in the Lie algebra $\frak g^{\Bbb C}$ is the subalgebra of all
holomorphic vector fields that admit an extension to the geometrical image
of the morphism. The action $\frak g^{\Bbb C}$ on the set of morphisms can
be exponentiated to the action of the semigroup
$[\Mantle(\Diff_+(S^1))]^{\#\frak A+\#\frak B}$, with $\#\frak A$ copies of the
semigroup $\Mantle(\Diff_+(S^1))$
acting from the right, and $\#\frak b$ from the left. Although the action of
the Lie algebra $\frak g^{\Bbb C}$ on the set of morphisms is not transitive,
it is possible to define the concept of a representation of the family of
morphisms as a continuous family of representations of the orbits of this
action.

A (projective) representation of the category $\Train(\Diff_+(S^1))$,
the train of the group $\Diff_+(S^1)$ ({\it (conformal) modular
functor\/}), is a set of representations $\pi_{\frak A\frak B}:\Mor(\frak
A,\frak B)\mapsto\BP(\Hom(V^{\otimes\#\frak A},V^{\otimes\#\frak B}))$ such
that
\roster
\item"--" the diagram
$$
\CD
\Mor(\frak A,\frak B)\times\Mor(\frak B,\frak C)@>>s>\Mor(\frak A,\frak C)\\
@VV\pi_{\frak A\frak B}\times\pi_{\frak B\frak C}V @VV\pi_{\frak A\frak C}V\\
\BP(\Hom(V^{\otimes\#\frak A},V^{\otimes\#\frak B})\times
\BP(\Hom(V^{\otimes\#\frak B},V^{\otimes\#\frak C}) @>>>
\BP(\Hom(V^{\otimes\#\frak A},V^{\otimes\#\frak C})),
\endCD
$$
where the lower arrow is the contraction operation, is commutative;
\item"--" if $\frak A=\frak A_1\sqcup\frak A_2$, $\frak B=\frak B_1\sqcup
\frak B_2$, $R=R_1\sqcup R_2$, $R_i\in\Mor(\frak A_i,\frak B_i)$, then
$\pi_{\frak A\frak B}(R)=\pi_{\frak A_1\frak B_1}(R_1)\times
\pi_{\frak A_2\frak B_2}(R_2)$.
\endroster

Every conformal modular functor corresponds to some representation of
poly\-nions, since polinions are a special case of morphisms in the
projective category $\Train(\Diff_+(S^1))$, the train of the group
$\Diff_+(S^1)$. In general, the converse is not true -- not every
representation of polynions can be extended to a conformal modular functor.
Indeed, the semigroup $\Voile(\Diff_+(S^1))$, the veil of the group
$\Diff_+(S^1)$, is identified with the semigroup of all endomorphisms of object
of cardinality 1 in the category $\Train(\Diff_+(S^1))$, and this is the
``topological'' obstruction to an extension of representations of polynions to
projective modular functors. Thus, among the structures of the theory of
represesentations in {\it the quantum conformal field theory\/}
[QPFT-operator algebras, representations of (conformal) polynions, (conformal)
modular functors] the last [the representations of the category
$\Train(\Diff_+(S^1))$, the train of the group $\Diff_+(S^1)$] form the
smallest class (under the condition that the representation spaces are sums of
Verma modules over the Lie algebra $\Cvir$), and the first [the QCFT-operator
algebras] form the largest class.

Note that each conformal modular functor is a projective modular functor
as well as each (projective) representation of (conformal) polynions
is a (projective) representation of the projective polynions.

\head Topic Five: $q_R$--conformal modular functors
\endhead

Let us now plunge the infinite dimensional picture of [1] into the
framework of two preceeding topics. However, a preliminary lemma is
necessary. In this topic the representation means the projective
representation.

\proclaim{Lemma}\newline
{\bf A.} The QPFT-operator algebra $\frak V$ generated by the
$\sltwo$--primary fields of non-negative integral spins $n$ (with the field
of spin 2 as $q_R$--conformal stress-energy tensor) [26] in the Verma module
$V_h$ of the extremal weight $h$ determines the representation of projective
polynions in the space $\frak V=\oplus_n V_n$.\newline
{\bf B.} The representation of projective polynions in the space $\frak V$
is naturally extended to the projective modular functor.
\endproclaim

The first statement of the lemma is almost evident. To prove the second
statement one should mark that the QPFT-operator algebra $\frak V$ may be
supplied with the operators $r_u(\varphi)$ of the multiplication from the right:
$r_u(\varphi)\psi=l_{-u}(\psi)\varphi$. Such operators generate the
QPFT-operator algebra and commute with the operators $l_u(\varphi)$ of the
initial QPFT-operator algebra $\frak V$. Both structures of QPFT-operator
algebras (of vertices $l_u(\varphi)$ and co-vertices [21] $r_u(\varphi)$)
supply the space $\frak V$ with the structure of the QPFT-operator
crossing-algebra, from which the projective modular functor may be restored
[21].

Let us put $\widetilde{\frak V}=\frak V\oplus V_h$ ($V_h$ -- unitarizable
Verma module). The space $\widetilde{\frak V}$ is supplied with the structure
of the QPFT-operator algebra, which is an abelian extension of $\frak V$.

\remark{Corollary} The QPFT-operator algebra $\widetilde{\frak V}$ determines
the representation of projective polynions in the space $\widetilde{\frak V}$,
which is an extension of the representation of projective polynions in
$\widetilde{\frak V}$.
\endremark

Let $\pi$ be a $\Cal K$--pseudorepresentations of the semigroup
$\Mantle(\Diff_+(S^1))$ in the space $V$ [1] (here $\Cal K$ denotes the class
of the compact operators in $V$). The class $\Cal K$ includes the class
$\operatorname{\Cal H\Cal S}$ of all Hilbert-Schmidt operators so the
$\operatorname{\Cal H\Cal S}$--pseu\-do\-rep\-re\-sen\-tations of [1] are
automatically $\Cal K$--pseudorepresentations.

The projective modular functor in $\frak V$ define a representation of the
projective semigroup $\Voile(\PSLTWO)$ in this space.

\proclaim{Theorem A}\newline
The representation of the projective semigroup
$\Voile(\PSLTWO)$ in the space $\frak V$ may be extended to the
$\Cal K$--pseudorepresentation of the semigroup $\Voile(\Diff_+(S^1))$
compatible with the $\Cal K$--representation of the semigroup
$\Mantle(\Diff_+(S^1))$.
\endproclaim

For any linear operator $A$
from $\Hom(V^{\otimes n};V)$ let us define the operators $F_{\bold x}^{(i)}(A)$
from $\End(V)$ ($1\le i\le n$, $\bold x\in V^{\otimes(n-1)}$), which are
received from $A$ by the substitution of $\bold x$ instead of all arguments
except the $i$-th argument. The linear operator $A$ from
$\Hom(V^{\otimes n};V)$ will be called {\it polycompact\/} if it is a linear
combination of the operators $A^{(i)}$ such that the operators
$F_{\bold x}^{(i)}(A)$ are compact for all $\bold x$. The (projective)
$\PK$--pseu\-do\-rep\-re\-sen\-tation of (conformal) polynions is a family
of representations $\pi_n$ of the homogeneous manifolds
$\Polynion_n(\Diff_+(S^1))$ up to the polycompact operators compatible with
$\pi$:
$$\pi_n:\Polynion_n(\Diff_+(S^1))\mapsto\BP(\Hom_{\Cal
B/\PK}(V^{\otimes(n+1)};V)),$$
where $\Hom_{\Cal B/\PK}$ is the quotient of the space of all bounded
operators by the subspace of all polycompact operators, such that the diagram
$$
\CD
\Polynion_n(\Diff_+(S^1))\times\Polynion_m(\Diff_+(S^1))@>>s>
\Polynion_{n+m}(\Diff_+(S^1))\\
@VV{\pi_n\times\pi_m}V                             @VV{\pi_{n+m}}V\\
\BP(\Hom_{\Cal B/\PK}(V^{\otimes(n+1)};V)\times\BP(\Hom_{\Cal
B/\PK}(V^{\otimes(m+1)};V)@>>>\BP(\Hom_{\Cal B/\PK}(V^{\otimes(n+m+1)};V))
\endCD
$$
where the lower arrow is the contraction operation, is commutative.

\proclaim{Theorem B}\newline
The representation of the projective polynions in the space $\frak V$
may be extended to the $\PK$--pseudorepresentation of the conformal
polynions compatible with the $\Cal K$--pseudorepresentation of the semigroup
$\Mantle(\Diff_+(S^1))$ in this space.
\endproclaim

\remark{Corollary}
The representation of the projective polynions in the space
$\widetilde{\frak V}$ of the corollary to lemma may be extended to the
$\PK$--pseudorepresentation of the conformal polynions compatible
with the $\Cal K$--pseudorepresentation of the semigroup
$\Mantle(\Diff_+(S^1))$ in $\widetilde{\frak V}$.
\endremark

The construction of the $\PK$--pseudorepresentations of polynions may be
reformulated for the antipolynions.

For any the linear operator $A$ from $\Hom(V^{\otimes n},V^{\otimes m})$
let us denote by $G^{(j)}_{\bold y}(A)$ ($1\le j\le m$, $\bold y\in
V^{\otimes(m-1)}$) the operator from $\Hom(V^{\otimes n},V)$, which is
received from $A$ by the pairing of the image of $A$ with $y$ by all
variables except the $j$-th one. Such linear operator $A$ will be called
{\it polycompact\/} if it is a linear combination of the operators $A^{(j)}$
such that $G^{(j)}_{\bold y}(A^{(j)})$ are polycompact operators uniformly
by $\bold y$ from any bounded subset of $V^{\otimes(m-1)}$.

The (projective) semi--$\PK$--pseudorepresentation of the category
$\Train(\Diff_+(S^1))$, the train of the group $\Diff_+(S^1)$, is a set of
representations $\pi_{\frak A\frak B}:\Mor(\frak A,\frak
B)\mapsto\BP(\Hom_{\Cal B/\PK}(V^{\otimes\#\frak A},V^{\otimes\#\frak
B}))$ up to the polycompact operators such that
\roster
\item"--" the set $\pi_{\frak A\frak B}$ being restricted to
$\Train(\PSLTWO)$ realizes the projective modular functor;
\item"--" the diagram
$$
\CD
\Mor(\frak A,\frak B)\times\Mor(\frak B,\frak C)@>>s>\Mor(\frak A,\frak C)\\
@VV\pi_{\frak A\frak B}\times\pi_{\frak B\frak C}V @VV\pi_{\frak A\frak C}V\\
\BP(\Hom_{\Cal B/\PK}(V^{\otimes\#\frak A},V^{\otimes\#\frak B})\times
\BP(\Hom_{\Cal B/\PK}(V^{\otimes\#\frak B},V^{\otimes\#\frak C}) @>>>
\BP(\Hom_{\Cal B/\PK}(V^{\otimes\#\frak A},V^{\otimes\#\frak C})),
\endCD
$$
where the lower arrow is the contraction operation, is commutative for any
two morphisms $f_1\in\Mor(\frak A,\frak B)$, $f_2\in\Mor(\frak B,\frak C)$ such
that $g(f_1\circ f_2)=g(f_1)+g(f_2)$, where $g(f)$ is the genus of the
geometric image of the morphism $f$;
\item"--" if $\frak A=\frak A_1\sqcup\frak A_2$, $\frak B=\frak B_1\sqcup
\frak B_2$, $R=R_1\sqcup R_2$, $R_i\in\Mor(\frak A_i,\frak B_i)$, then
$\pi_{\frak A\frak B}(R)=\pi_{\frak A_1\frak B_1}(R_1)\times
\pi_{\frak A_2\frak B_2}(R_2)$.
\endroster

\proclaim{Theorem C}\newline
The projective modular functor in the space $\frak V$ may be extended to the
semi--$\PK$--pseu\-do\-rep\-re\-sen\-tation of the category
$\Train(\Diff_+(S^1))$ compatible with the
$\PK$--pseu\-do\-rep\-re\-sen\-tations of the conformal polynions and
antipolynions as well as with the $\Cal K$--pseudorepresentation of the
semigroup $\Voile(\Diff_+(S^1))$.
\endproclaim

The semi--$\PK$--pseudorepresentation of the category $\Train(\Diff_+(S^1))$
will be called the {\it $q_R$--conformal modular functor}. It may be considered
as the ``Berezin quantization'' of the conformal modular functor in view of
the original construction of the $q_R$--conformal symmetries and QPFT--operator
algebra $\frak V$ from the Lobachevski{\v\i} $C^*$-algebra, the Berezin
quantization of the Lobachevski{\v\i} plane [26] (see also [4]); here
$q_R=\tfrac1{2h-1}$ is the quantization parameter.

\remark\nofrills{Problems:}
\roster
\item"--" To investigate the asymptotic behaviour of the
$\PK$--pseudorepresentations of conformal polynions and $q_R$--conformal
modular functors in the space $\frak V$ if as $h$ tends to $\infty$
(or $q_R$ tends to zero).
\item"--" To explore the possible relations between $q_R$--conformal
modular functors and quantizations of Riemann surfaces in sense of
S.Klimek and A.Les\-ni\-ew\-s\-ki [50].
\endroster
\endremark

\Refs
\roster
\item" [1]" Juriev D., On the infinite-dimensional hidden symmetries. I.
Infinite dimensional geometry of $q_R$-conformal symmetries. E-print:
funct-an/9612004.
\item" [2]" Juriev D., Topics in hidden symmetries. I-V. E-prints:
hep-th/9405050, q-alg/9610026, q-alg/9611003, q-alg/9611019, funct-an/9611003.
\item" [3]" Juriev D.V., An excursus into the inverse problem of representation
theory [in Russian]. Report RCMPI-95/04 [e-version: mp\_arc/96-477].
\item" [4]" Juriev D.V., Octonions and binocular mobilevision [in Russian].
Fundam.Prikl.Matem. 3 (1997), to appear; Visualizing 2D quantum field theory:
geometry and in\-for\-matics of mobilevision. Report RCMPI/96-02 [{\sl draft}
e-version: hep-th/9401067+9403137]; Be\-lav\-kin-Ko\-lo\-kol\-tsov watch-dog
effects in interactively controlled stochactic dynamical videosystems [in
Russian]. Teor.Matem.Fiz. 106(2) (1996) 333-352 [English transl.: Theor.Math.
Phys. 106 (1996) 276-290]; On the description of a class of physical
interactive information systems [in Russian]. Report RCMPI/96-05 [e-version:
mp\_arc/96-459].
\item"[5]" Juriev D., Topics in hidden algebraic structures and infinite
dimensional dynamical symmetries of controlled systems, in preparation.
\item"[6]" Segal G., Definitions of conformal field theory, Preprint MPI, 1988.
\item"[7]" Gawedzki K., Conformal field theory. Sem.Bourbaki, 1988, ex.704.
\item"[8]" Moore G., Seiberg N., Classical and quantum conformal field
theory. Commun.Math. Phys. 123 (1989) 177-254.
\item"[9]" Neretin Yu.A., A complex semigroup that contains the group of
diffeomorphisms of the circle [in Russian]. Funkts.anal.i ego prilozh. 21(2)
(1987) 82-83; Holomorphic extensions of representations of the
group of diffeomorphisms of the circle [in Russian]. Matem.Sbornik 180(5)
(1989) 635-657; Infinite-dimensional groups, their mantles, trains and
representations. Adv.Soviet Math. 2 (1991) 103-171; Ca\-te\-go\-ries of
symmetries and infinite-dimensional groups. London Math.Soc.Mono. 16,
Cla\-ren\-don Press, 1996.
\item"[10]" Witten E., Quantum field theory, grassmannians and algebraic curves.
Commun.Math. Phys. 113 (1988) 529-600; Physics and geometry. Rep. 100th
Anniversary A.M.S., 1988.
\item"[11]" Arbarello E., De Concini C., Kac V., Procesi C., Moduli spaces
of curves and representation theory. Commun.Math.Phys. 117 (1988) 1-36;
Alvarez-Gaume L., Gomez G., Moore G., Vafa C., Strings in operator formalism.
Nucl.Phys. B303 (1988) 455-521.
\item"[12]" Zhelobenko D.P., Representations of the reductive Lie algebras
[in Russian]. Moscow, Nauka, 1994.
\item"[13]" Sklyanin E.K., Takhtajan L.A., Faddeev L.D., Quantum inverse
problem method [in Russian]. Teor.Matem.Fiz. 40(2) (1979) 194-220; Sklyanin
E.K., Quantum method of inverse scattering problem [in Russian].
Zap.Nauch.Semin.LOMI 95 (1980) 55-128; Kulish P.P., Sklyanin E.K., Quantum
spectral transform method. Recent developments. Lect.Notes Phys. 151 (1982)
61-119; Izergin A.G., Korepin V.E., Quantum inverse problem method [in
Russian]. Phys.Elem.Particles \& Atom.Nuclei 13(3) (1982) 501-541.
\item"[14]" Karasev M.V., Maslov V.P., Naza{\v\i}kinski{\v\i} V.E., Algebras
with general commutation relations [in Russian]. Current Probl.Math.,
Modern Achievements. V.13. Moscow, VINITI, 1979.
\item"[15]" Karasev M.V., Maslov V.P., Nonlinear Poisson brackets. Geometry
and quantization. Amer.Math.Soc., Providence, RI, 1993.
\item"[16]" Juriev D.V., Quantum projective field theory: quantum-field
analogs of Euler-Arnold equations in projective $G$--hypermultiplets [in
Russian]. Teor.Matem.Fiz. 98(2) (1994) 220-240 [English transl.:
Theor.Math.Phys. 98 (1994) 147-161].
\item"[17]" Gelfand I.M., Gindikin S.G., Complex manifolds, whose frames are
semisimple real Lie groups, and analytic discrete series of representations
[in Russian]. Funkts.anal.i ego pri\-lozh. 11(1) (1977) 19-27; Vinberg E.B.,
Invariant convex cones and orderings in Lie groups [in Russian]. Funkts.anal.i
ego prilozh. 14(1) (1980) 1-13; Olshanski{\v\i} G.I., Invariant cones in Lie
algebras, Lie semigroups and holomorphic discrete series [in Russian].
Funkts.anal.i ego prilozh. 15(4) (1981) 53-66.
\item"[18]" Hilgert J., Hofmann K.H., Lawson J., Lie groups, convex cones and
semigroups. Clarendon Press, Oxford, 1989.
\item"[19]" Gelfand I.M., Minlos R.A., Shapiro Z.Ya., Representations of the
rotation and Lorentz groups and their applications. Pergamon Press, Oxford,
1963; Gelfand I.M., Graev M.I., Vilenkin N.Ya., Integral geometry and
representation theory. Acad.Press, New York, 1966; Vilenkin N.Ya., Special
functions and the theory of group representations, A.M.S. Transl.Math.Mono.
22, Providence, RI, 1968; Shelepin L.A., Clebsch-Gordan coefficient calculus
and its physical applications [in Russian]. Trans.Phys.Inst. USSR Acad.Sci. 70
(1973) 3-119; Klimyk A.U., Matrix elements and Clebsch-Gordan coefficients of
group representations [in Russian]. Naukova Dumka, Kiev, 1979.
\item"[20]" McKean H.P., Trubowitz E., Hill's operator and hyperelliptic
function theory in the presence of infinitely many branch points.
Commun.Pure Appl.Math. 29(2) (1976) 143-226.
\item"[21]" Juriev D., Infinite dimensional geometry and quantum field theory
of strings. II. Infinite-dimensional noncommutative geometry of
self-interacting string field. Russian J.Math. Phys. 4(3) (1996).
\item"[22]" Kirillov A.A., Elements of the representation theory. Springer,
1976.
\item"[23]" Barut A., Raczka R., Theory of group representations and
applications, Warszawa, 1977.
\item"[24]" Juriev D.V., Classification of the vertex operators in
two-dimensional $\sltwo$-invariant quantum field theory [in Russian].
Teor.Matem.Fiz. 86(3) (1991) 338-343; The explicit form of the vertex
operator fields in two-dimensional quantum $\sltwo$-invariant field theory.
Lett.Math.Phys. 22 (1991) 141-144.
\item"[25]" Bychkov S.A., Juriev D.V., Three algebraic structures of the
quantum projective ($\sltwo$-invariant) field theory [in Russian].
Teor.Matem.Fiz. 97(3) (1993) 336-347 [English transl.: Theor.Math.Phys. 97
(1993) 1333-1339].
\item"[26]" Juriev D.V., Complex projective geometry and quantum projective
field theory [in Russian]. Teor.Matem.Fiz. 101(3) (1994) 331-348 [English
transl.: Theor.Math.Phys. 101 (1994) 1387-1403].
\item"[27]" Kirillov A.A., An introduction to the representation theory and
noncommutative harmonic analysis [in Russian]. Current Probl.Math.,
Fundam.Directions, V.22, Moscow, VINITI, 1988.
\item"[28]" Juriev D.V., QPFT-operator algebras and commutative exterior
differential calculus [in Russian]. Teor.Matem.Fiz. 93(1) (1992) 32-38
[English transl.: Theor.Math.Phys. 93 (1992) 1101-1105].
\item"[29]" Juriev D.V., Quantum conformal field theory as infinite-dimensional
noncommutative geometry [in Russian]. Uspekhi Matem.Nauk 46(4) (1991) 115-138
[English transl.: Russian Math.Surveys 46(4) (1991) 135-163].
\item"[30]" Krichever I.M., Novikov S.P., Virasoro-type algebras, Riemann
surfaces and structures of the soliton theory [in Russian]. Funkts.analiz i
ego prilozh. 21(2) (1987) 46-63; Virasoro-type algebras, Riemann surfaces
and string in the Minkowsky space [in Russian]. Funkts. anal.i ego prilozh.
21(4) (1987) 47-61; Virasoro-type algebras, the momentum-energy tensor and
operator expansions on the Riemann surfaces [in Russian]. Funkts.anal.i ego
prilozh. 23(1) (1989) 1-14; Virasoro-Gelfand-Fuchs algebras, Riemann surfaces,
operator's theory of closed strings. J.Geom.Phys. 5 (1988) 631-661 [reprinted
in ``Geometry and physics. Essays in honour of I.M.Gelfand'', Eds. S.Gindikin
and I.M.Singer, Pitagora Editrice, Bologna and Elsevier Sci.Publ., Amsterdam,
1991].
\item"[31]" Sklyanin E.K., On algebraic structures related to the Yang-Baxter
equation [in Russian]. Funkts. analiz i ego prilozh. 16(4) (1982) 27-34; 17(4)
(1983) 34-48; On an algebra generated by the quadratic relations [in Russian].
Uspekhi Matem. Nauk 40(2) (1985) 214.
\item"[32]" Sabinin L.V., On the nonlinear geometric algebra [in Russian]. In
``Webs and quasigroups''. Kalinin [Tver'], 1988, pp.32-37; Methods of
nonassociative algebra in differential geometry [in Russian]. Suppl.to the
Russian transl.of Kobayashi S., Nomizu K., Foundations of differential geometry.
V.1. pp.293-339, Moscow, Nauka, 1982; Differential geometry and quasigroups [in
Russian]. In ``Current problems of geometry. To the 60th Anniversary of
Acad.Yu.G.Reshetnyak''. Trans.Inst.Math. Siberian Branch Soviet Acad.Sci.,
Novosibirsk, 1984, v.14, pp.208-221; Mikheev P.O., Sabinin L.V., Quasigroups
and differential geometry. In ``Quasigroups and loops. Theory and
applications''. Berlin, Heldermann Verlag, 1990, P.357-430.
\item"[33]" Juriev D.V., Topics in isotopic pairs and their representations
[in Russian]. Teor.Matem. Fiz. 105(1) (1995) 18-28 [English transl.:
Theor.Math.Phys. 105 (1995) 1201-1209]; Topics in isotopic pairs and their
representations. II. A general supercase [in Russian]. Teor.Ma\-tem.Fiz.
111(1) (1997), to appear.
\item"[34]" Reshetikhin N.Yu., Takhtajan L.A., Faddeev L.D., Quantization of
Lie groups and Lie algebras [in Russian]. Algebra i analiz 1(1) (1989) 178-206
[English transl.: St.Petersburg Math.J. 1 (1990) 193-225].
\item"[35]" Ro\v cek M., Representation theory of the nonlinear $\SU(2)$
algebra. Phys. Lett.B 255 (1991) 554-557.
\item"[36]" Olshanski{\v\i} G.I., Method of holomorphic extensions in representation
theory of infinite-dimensional classical groups [in Russian]. Funkts.anal.i ego
prilozh. 22(4) (1988) 23-37; On semigroups related to infinite-dimensional
groups. Adv.Soviet Math. 2 (1991) 67-101.
\item"[37]" Gelfand I.M., Fuchs D.B., The cohomology of the Lie algebra of
vector fields on a circle [in Russian]. Funkts.anal.i ego prilozh. 2(4) (1968)
92-93.
\item"[38]" Fuchs D.B., Cohomology of infinite dimensional Lie algebras
[in Russian]. Moscow, Nauka, 1984.
\item"[39]" Virasoro M.A., Subsidiary conditions and ghosts in dual resonance
models. Phys.Rev.D. 1 (1970) 2933-2936.
\item"[40]" Bott R., The characteristic classes of groups of diffeomorphisms.
Enseign.Math. 23 (1977) 209-220.
\item"[41]" Juriev D., Infinite dimensional geometry and quantum field theory
of strings. I. Infinite-dimensional geometry of second quantized free string.
Alg.Groups Geom. 11 (1994) 145-179.
\item"[42]" Segal G., Unitary representations of some infinite dimensional
groups. Commun.Math. Phys. 80 (1991) 301-342; Kirillov A.A.,
Infinite-dimensional Lie groups, their orbits, in\-variants and representations.
Lect.Notes Math. 970, Springer, 1982; K\"ahler structure on K-orbits of the
group of diffeomorphisms of a circle [in Russian]. Funkts.anal.i ego prilozh.
21(2) (1987) 42-45; Kirillov A.A., Juriev D.V., K\"ahler geometry of the
infinite dimensional space $M=\Diff_+(S^1)/\Rot(S^1)$ [in Russian].
Funkts.anal.i ego prilozh. 21(4) (1987) 35-46; Representations of the
Virasoro algebra by the orbit method. J.Geom.Phys. 5 (1988) 351-363
[reprinted in ``Geometry and Physics. Essays in honour of I.M.Gelfand''.
Eds.S.G.Gindikin and I.M.Singer, Pitagora Editrice, Bologna and Elsevier
Sci.Publ., Amsterdam, 1991].
\item"[43]" Juriev D., The vocabulary of geometry and harmonic analysis on the
infinite-dimensional manifold $\Diff_+(S^1)/S^1$. Adv.Soviet Math. 2 (1991)
233-247; A model of the Verma modules over the Virasoro algebra [in Russian].
Algebra i analiz 2(2) (1990) 209-226; Non-Euclidean geometry of mirrors and
prequantization on the homogeneous K\"ahler manifold $M=\Diff_+(S^1)/\Rot(S^1)$
[in Russian]. Uspekhi Matem.Nauk 43(2) (1988) 159-160;
Infinite-dimensional geometry of the universal deformation of the complex
disk. Russian J.Math.Phys. 2(1) (1994) 111-121.
\item"[44]" Kac V.G., Infinite dimensional Lie algebras. Cambridge, Cambridge
Univ. Press, 1990.
\item"[45]" Mack G., Introduction to conformal invariant quantum field theory
in two and more dimensions. In ``Nonperturbative quantum field theory''. Eds.
G.t'Hooft et al., Plenum, New York, 1988;
Hadjivanov L.K., Existence of primary fields as a generalization
of the Luscher-Mack theorem. J.Math.Phys. 34 (1993) 441-453.
\item"[46]" Juriev D.V., The algebra $\Vrt(\Cvir;c)$ of vertex operators for
the Virasoro algebra [in Russian]. Algebra i analiz. 3(3) (1991) 197-205;
Representation of operator algebras of the quantum conformal field theory.
J.Math.Phys. 33 (1992) 492-496.
\item"[47]" Frenkel I., Lepowsky J., Meurman A., Vertex operator algebras
and the Monster. Acad. Press, 1988; Dong C., Lepowsky J., Generalized
vertex algebras and relative vertex operators. Birkh\"auser, 1993;
Kac V.G., Vertex algebras for Beginners. Amer.Math.Soc., Providence, RI,
1996.
\item"[48]" Juriev D., Local conformal field algebras. Commun.Math.Phys. 138
(1991) 569-581, (E) 146 (1992) 427; On the structure of L-algebra $L(\Cvir)$.
J.Funct.Anal. 101 (1991) 1-9.
\item"[49]" Juriev D., Projective reduction in quantum conformal field theory.
J.Math.Phys. 32 (1991) 2034-2038.
\item"[50]" Klimek S., Lesniewski A., Quantum Riemann surfaces. I-III.
Commun.Math.Phys. 146 (1992) 103-122, Lett.Math.Phys. 24 (1992) 125-139,
32 (1994) 45-61.
\endroster
\endRefs
\enddocument